\documentclass{amsproc}
\usepackage{amsmath,amsthm,amsfonts}

\newtheorem{theorem}{Theorem}[section]
\theoremstyle{definition}
\newtheorem{definition}[theorem]{Definition}
\newtheorem{corollary}[theorem]{Corollary}
\newtheorem{remark}[theorem]{Remark}
\newtheorem{lemma}[theorem]{Lemma}
\numberwithin{equation}{section}
\def\dot{\bf\mathaccent"705F }
\newcommand{\bkRp}{{\bf \dot{{\mathbb{R}}}}}
\newcommand{\dist}{\mbox{\rm dist}}

\begin{document}

\title[Wiener-Hopf plus Hankel operators]{Wiener-Hopf plus Hankel operators on the real line\\
with unitary and sectorial symbols}

\author{G. Bogveradze}
\address{Department of Mathematics, University of Aveiro,
3810-193 Aveiro, Portugal} \email{giorgi@mat.ua.pt}
\thanks{G. Bogveradze is sponsored by a grant from the Research Unit {\it Mathematics and Applications}
through {\em Funda\c{c}{\~a}o para a Ci{\^e}ncia e a Tecnologia} (Portugal).}

\author{L.$\!$ P.~Castro}
\address{Department of Mathematics, University of Aveiro,
3810-193 Aveiro, Portugal} \email{lcastro@mat.ua.pt}
\urladdr{http://www.mat.ua.pt/lcastro/}



\maketitle

\begin{abstract}
Wiener-Hopf plus Hankel operators acting between Lebesgue spa\-ces
on the real line are studied in view of their invertibility, one
sided-in\-ver\-ti\-bi\-li\-ty, Fredholm, and semi-Fredholm
properties. This is done in two different cases: (i) when the
Fourier symbols of the operators are unitary functions, and (ii)
when the Fourier symbols are related with sectorial elements
appearing in factorizations of functions originated by the Fourier
symbols of the operators. The obtained result for the case (i) may
be viewed as a Douglas-Sarason type theorem for Wiener-Hopf plus
Hankel operators.
\end{abstract}

\section{Introduction}

Algebraic sums of Wiener-Hopf and Hankel operators have been
receiving increasing attention in the last years
(cf.~\cite{BaTo04}, \cite{CaSp05pm}, \cite{CaSpTe04mn},
\cite{CoKrTe03}, \cite{Eh04}--\cite{NoCa05b}, \cite{RoSi90},
\cite{Te89}). A great part of the interest is directly originated
by concrete mathematical-physics applications where this kind of
operators appear. This is the case of problems in wave diffraction
phenomena which are modelled by boundary-transmission value
problems that can be equivalently translated into systems of
integral equations characterized by such kind of operators (see,
e.g., \cite{CaSpTe03jie}--\cite{CaSpTe06ieot}).

It is therefore easy to understand that a fundamental goal in those studies
consists in obtaining conditions that lead to the {\it regularity properties}
of the operators (i.e., their invertibility, Fredholm property, and other
properties that are directly dependent on the kernel and image of the
operators). For some classes of the so-called {\it Fourier symbols of the
operators} their regularity properties are already known. This is the case of
the continuous or, even, piecewise continuous functions on the compactificated
real line (cf.~the above references). Despite those advances, for some other
classes of Fourier symbols a complete description of the regularity properties
is still missing.

In the present paper, starting simply with essential bounded Fourier symbols,
we look for conditions that may ensure the above pointed properties for
Wiener-Hopf plus Hankel operators having that Fourier symbols.  All the
detailed definitions are given in the next section.

\section{The Wiener-Hopf plus Hankel operators in study}

We will consider Wiener-Hopf plus Hankel operators of the form
\begin{eqnarray}
\label{Op:WHH1}
&&WH_{\varphi}=W_{\varphi}+H_{\varphi} :
L_{+}^{2}(\mathbb{R}) \rightarrow L^{2}(\mathbb{R}_{+})\;,
\end{eqnarray}
with
$W_{\varphi}$ and $H_{\varphi}$ being Wiener-Hopf and Hankel
operators defined by
\begin{eqnarray}\label{Op:WHH2}
&&W_{\varphi}=r_{+}\mathcal{F}^{-1}\varphi\cdot\mathcal{F} :
L_{+}^{2}(\mathbb{R}) \rightarrow L^{2}(\mathbb{R}_{+})\\
&&H_{\varphi}=r_{+}\mathcal{F}^{-1}\varphi\cdot\mathcal{F}J :
L_{+}^{2}(\mathbb{R}) \rightarrow
L^{2}(\mathbb{R}_{+})\;,\label{Op:WHH3}
\end{eqnarray}
respectively. As usual, $L^2(\mathbb{R})$ and $L^2(\mathbb{R}_+)$
denote the Hilbert spaces of complex-valued Lebesgue measurable
functions $\varphi$, for which ${|\varphi|}^2$ is integrable on
$\mathbb{R}$ and $\mathbb{R}_+$, respectively. In
(\ref{Op:WHH1})--(\ref{Op:WHH3}), $L_{+}^{2}(\mathbb{R})$ denotes
the subspace of $L^{2}(\mathbb{R})$
 formed by all functions supported in the closure of
$\;\mathbb{R}_{+}=(0,+\infty)$, the operator $r_{+}$ performs the
restriction from $L^{2}(\mathbb{R})$ into $L^{2}(\mathbb{R_{+}})$,
$\mathcal{F}$ denotes the Fourier transformation, $J$ is the
reflection operator given by the rule
$J\varphi(x)=\widetilde\varphi(x)=\varphi(-x)$, $x\in
\mathbb{\mathbb{R}}$, and $\varphi\in L^\infty(\mathbb{R})$ is the
so-called Fourier symbol.

Let $\mathbb{C_-}=\{z \in \mathbb{C}: {\Im}m\,z<0\}$ and $\mathbb{C_+}=\{z \in
\mathbb{C}: {\Im}m\,z>0\}$. As usual, let us denote by
$H^{\infty}(\mathbb{C_\pm})$ the set of all bounded and analytic functions in
$\mathbb{C_\pm}$. {\em Fatou's Theorem} ensures that functions in
$H^{\infty}(\mathbb{C_\pm})$ have non-tangential limits on $\mathbb{R}$ almost
everywhere. Thus, let $H_{\pm}^{\infty}(\mathbb{R})$ be the set of all
functions in $L^{\infty}(\mathbb{R})$ that are non-tangential limits of
elements in $H^{\infty}(\mathbb{C_\pm})$. Moreover, it is well known that
$H_{\pm}^{\infty}(\mathbb{R})$ are closed subalgebras of
$L^\infty(\mathbb{R})$.

We will concentrate on obtaining regularity properties of this kind of
Wiener-Hopf plus Hankel operators (\ref{Op:WHH1}) in the case where the Fourier
symbol $\varphi$ is  an unitary-valued function or a function somehow
associated to a sectorial element through a certain factorization (and in these
two cases we will make use of the just presented $H_{\pm}^{\infty} :=
H_{\pm}^{\infty}(\mathbb{R})$ spaces). We recall that $\phi\in
L^\infty(\mathbb{R})$ is called {\em unitary} if $|\phi| = 1$, and is called
{\em sectorial} if the essential range of $\phi$ is contained in an open
half-plane whose boundary passes through the origin.

As mentioned, within the context of the present paper we are using
the name of {\em regularity properties} of a linear operator
acting between Hilbert spaces for those proprieties that arise
from a direct influence of the kernel and the image of those
operators. In more detail, let $T:X\rightarrow Y$ be a bounded
linear operator acting between Hilbert spaces and, if
$\text{Im}\,T$ is closed (i.e. $T$ is {\em normally solvable}),
let us consider the cokernel of $T$ to be defined by the quotient
$\text{Coker}\,T=Y/ \text{Im}\,T$. Then, a normally solvable
operator $T$ is said to be {\em right-Fredholm} if $\dim
\text{Coker }T$ is finite, {\em left-Fredholm} if $\dim
\text{Ker}\,T$ is finite, and {\em Fredholm} if both $\dim
\text{Ker}\,T$ and $\dim \text{Coker}\,T$ are finite.
Additionally, we say that $T$ is {\em left-invertible} or {\em
right-invertible} if there exist $T^-_l :Y \rightarrow X$ or
$T^-_r :Y \rightarrow X$ such that $T^-_lT = I_X$ or $TT^-_r=
I_Y$, respectively. As usual, in the case when both $T^-_l$ and
$T^-_r$ exist and additionally $T^-_l = T^-_r$, the operator $T$
is said to be {\em invertible} (or {\em both sided-invertible}).
Alternatively, it can be shown that $T$ is left-invertible if and
only if $T$ is injective and normally solvable. In the same way,
$T$ is right-invertible if and only if $T$ is normally solvable
and surjective.

\section{A Douglas--Sarason type theorem}

For a Banach algebra $B$, we are going to denote by $\mathcal{G}B$ the group of
all invertible elements in $B$.

We point out that if $\phi\in\mathcal{G}L^\infty(\mathbb{R})$,
then there is a function $h \in \mathcal{G}H^\infty_+(\mathbb{R})$
such that $|\phi|=|h|$ almost everywhere in $\mathbb{R}$. This is
only a particular case of the {\em Pousson--Rabindranathan
Theorem}~\cite{Po68, Ra69} (cf.~Theorem~\ref{theorem PR} below)
but give a hint to the strategy behind the next results where
distances between certain classes and the Fourier symbols of the
operators are considered.

In the Wiener-Hopf operators case there is a well known theorem --
due to Douglas and Sarason -- about the Fredholm theory of this
kind of operators and the distances between the Fourier symbol to
certain spaces. More precisely, the theorem may be written in the
following form.

\begin{theorem}[Douglas and Sarason~\cite{DoSa}] \label{theorem1}
If $\varphi\in L^{\infty}(\mathbb{\mathbb{R}})$ is unitary, then:
\begin{itemize}
\item[(a)] $W_{\varphi}$ is invertible if and only if
$\;\dist(\varphi,\mathcal{G}H_{+}^{\infty})<1$ or if and only if
$\;\dist(\varphi,\mathcal{G}H^{\infty}_{-})<1\,;$\\[-3.003mm]
 \item[(b)]
$W_{\varphi}$ is left-invertible  if and only if
$\;\dist(\varphi,H^{\infty}_{+})<1\,;$\\[-3.003mm] \item[(b$^\prime$)]
$W_{\varphi}$ is right-invertible if and only if
$\;\dist(\varphi,H^{\infty}_{-})<1\,;$\\[-3.003mm] \item[(c)]$W_{\varphi}$ is
Fredholm if and only if
$\;\dist(\varphi,\mathcal{G}[C(\bkRp)+H^{\infty}_{+}])<1$ or if
and only if
$\;\dist(\varphi,\mathcal{G}[C(\bkRp)+H^{\infty}_{-}])<1\,;$\\[-3.003mm]
\item[(d)]$W_{\varphi}$ is left-Fredholm if and only if
$\;\dist(\varphi,C(\bkRp)+H^{\infty}_{+})<1\,;$\\[-3.003mm]
\item[(d$^\prime$)]$W_{\varphi}$ is right-Fredholm if and only if
$\;\dist(\varphi,C(\bkRp)+H^{\infty}_{-})<1$.
\end{itemize}
\end{theorem}

As usual, $C(\bkRp)$ denotes the space of all bounded continuous
(complex-valued) functions on $\mathbb{R}$ for which both limits at $\pm\infty$
exist and coincide.

The last theorem was the starting motivation for obtaining such kind of result
for our Wiener-Hopf plus Hankel operators. It is clear that adding an Hankel
operator to the above Wiener-Hopf operator will influence several changes in
the regularity properties of the resulting operator.

To consider such situation, we will use the notion of
$\Delta$-relation after extension introduced in \cite{CaSp98} for
bounded linear operators acting between Hilbert spaces, e.g. $T :
X_{1} \rightarrow X_{2}$ and $S : Y_{1} \rightarrow Y_{2}$. We say
that $T$ {\em is $\Delta$-related after extension to} $S$ (and use
the abbreviation $T\stackrel{*}{\Delta}S$) if there is an
auxiliary bounded linear operator acting between Hilbert spaces
$T_{\Delta} : X_{1\Delta} \rightarrow X_{2\Delta}$, and bounded
invertible operators $E$ and $F$ such that
\begin{equation}
\label{eq:delta}
\left[
  \begin{array}{cc}
    T & 0 \\
    0 & T_{\Delta} \\
  \end{array}
\right]= E \left[\begin{array}{cc} S & 0 \\
 0 & I_{Z}
 \end{array}
 \right]F\,,
 \end{equation}
 where $Z$ is an additional Hilbert space and $I_{Z}$ represents the
identity operator in $Z$. In the particular case where $T_{\Delta}
= I_{X_{1\Delta}} : X_{1\Delta} \rightarrow
X_{2\Delta}=X_{1\Delta}$ is the identity operator, we say the $T$
and $S$ {\em are equivalent after extension operators}. In the
simplest case of $T = E\, S\, F$ (for some boundedly invertible
operators $E$ and $F$) we say that $T$ and $S$ are {\em equivalent
operators}.

\begin{theorem}
Let $\varphi$ be unitary in $L^{\infty}(\mathbb{\mathbb{R}})$.
\begin{itemize} \item[(a)] $WH_{\varphi}$ and $W_{\varphi}-H_{\varphi}$ are invertible if
and only if
$$\dist(\varphi\widetilde{\varphi^{-1}},\mathcal{G}H^{\infty}_{+})<1$$ or if and
only if
$\;\dist(\varphi\widetilde{\varphi^{-1}},\mathcal{G}H^{\infty}_{-})<1\,.$\\[-2.003mm]
\item[(b)] $WH_{\varphi}$ and $W_{\varphi}-H_{\varphi}$ are
left-invertible if and only if
$$\dist(\varphi\widetilde{\varphi^{-1}}, H^{\infty}_{+})<1\,.$$
\item[(b$^\prime$)] $WH_{\varphi}$ and $W_{\varphi}-H_{\varphi}$
are right-invertible if and only if
$$\dist(\varphi\widetilde{\varphi^{-1}}, H^{\infty}_{-})<1\,.$$
\item[(c)] $WH_{\varphi}$ and $W_{\varphi}-H_{\varphi}$ are
Fredholm if and only if
$$\dist(\varphi\widetilde{\varphi^{-1}},\mathcal{G}[C(\bkRp)+H^{\infty}_{+}])<1\,,$$
or if and only if
$\dist(\varphi\widetilde{\varphi^{-1}},\mathcal{G}[C(\bkRp)+H^{\infty}_{-}])<1\,.$\\[-2.003mm]
\item[(d)]$WH_{\varphi}$ and $W_{\varphi}-H_{\varphi}$ are
left-Fredholm if and only if
$$\dist(\varphi\widetilde{\varphi^{-1}},C(\bkRp)+H^{\infty}_{+})<1\,.$$
\item[(d$^\prime$)] $WH_{\varphi}$ and $W_{\varphi}-H_{\varphi}$
are right-Fredholm if and only if
$$\dist(\varphi\widetilde{\varphi^{-1}},C(\bkRp)+H^{\infty}_{-})<1\,.$$
\end{itemize}
\end{theorem}

\begin{proof}
Our proof is based on the notion of $\Delta$-relation after
extension, which was recalled above for bounded linear operators
acting between Hilbert spaces. From \cite[Example 1.7]{CaSp98} we
can derive that $T = WH_{\varphi}$ is ${\Delta}$-related after
extension to the Wiener-Hopf operator $S =
W_{\varphi\widetilde{\varphi^{-1}}}$,
 cf.~also \cite[Lemma~2.1]{NoCa05} or \cite[Lemma~1]{NoCa05b} where
this relation is given in an explicit way for $T_{\Delta} =
W_\varphi - H_\varphi$ in (\ref{eq:delta}) (or for an operator
$T_{\Delta}$ equivalent to $W_\varphi - H_\varphi$). Thus, we are
now going to analyze the Fourier symbol
$\varphi\widetilde{\varphi^{-1}}$.

Let $\varphi^{*}$ stands for the complex conjugate of $\varphi$.
Since $\varphi$ is unitary, then
$\varphi\varphi^{*}=\varphi^{*}\varphi=1$, and whence
$\varphi^{-1}=\varphi^{*}$. From here we also obtain
${\widetilde{\varphi}^{-1}}={\widetilde{\varphi}^{*}}$. Observing
that
$$\widetilde{[\varphi\varphi^{*}]}=\widetilde{[\varphi^{*}\varphi]}=1\,,$$
we reach to
$\widetilde{\varphi}\widetilde{\varphi^{*}}=\widetilde{\varphi^{*}}\widetilde{\varphi}=1$.
This means that $\widetilde{\varphi}$ is unitary and so it is
$\varphi\widetilde{\varphi^{-1}}$. We can now apply Theorem~\ref{theorem1} to
the operator $W_{\varphi\widetilde{\varphi^{-1}}}$ and obtain all the above
stated conditions in terms of distances. Now, the result follows if we
interpret (\ref{eq:delta}) as an equivalence after extension relation between
$\mbox{diag}[T, T_{\Delta}]= \mbox{diag}[WH_{\varphi}, W_\varphi - H_\varphi]$
and $S= W_{\varphi\widetilde{\varphi^{-1}}}$.
\end{proof}

\begin{remark}{\rm
We would like to point out that in the above reasoning (of the
$\Delta$-relation after extension) it is not possible to exclude the
factor $\widetilde\varphi^{-1}$ for obtaining a ``pure" (direct)
result depending only on $\varphi$ (without extra factors). Consider
for instance the case (b). From the condition
$\,\dist(\varphi\widetilde{\varphi^{-1}},H^{\infty}_{+})<1\;$ it
follows that
$$\inf_{h\in H_{+}^{\infty}}{\|\varphi\widetilde{\varphi^{-1}}-h\|}_{L^{\infty}\!(\mathbb{\mathbb{R}})}<1\,,$$
which is equivalent to
$$\inf_{h\in H_{+}^{\infty}}{\|\varphi-h\widetilde{\varphi}\|}_{L^{\infty}\!(\mathbb{\mathbb{R}})}<1\;,$$
but this is the distance from $\varphi$ to the space
$\{h\widetilde{\varphi} : h \in H^{\infty}_{+} \;\mbox{and}\;
\varphi \;\mbox{is unitary}\}$. Note that this space is obviously
different from $H_{+}^{\infty}.$ In fact, if we take $\varphi\in
H_{+}^{\infty},$ with $\varphi(x)=e^{i\lambda x}$ and $\lambda>0$
$(x\in\mathbb{R}),$ then choosing $h(x)= e^{i\mu x}$ with
$0<\mu<\lambda$ we have $h\widetilde{\varphi}\in H_{-}^{\infty}.$
This also means that the space is directly dependent on the
structure of $\widetilde\varphi$ (and therefore on
$\widetilde{\varphi^{-1}}$).}
\end{remark}

\section{Sectorial Symbols}

In the present section we will work with the so-called {\it sectorial Fourier
symbols}.

\begin{definition}{\rm
A function $f\in L^{\infty}(\mathbb{R})$ is said to be {\it
sectorial} if there exist a real number $\varepsilon>0$ and a
$c\in\mathbb{C}$ of modulus $1$ such that
$${\Re}e (cf(x))\geq\varepsilon\;,$$
almost everywhere on  $\mathbb{R}$.}
\end{definition}

We will denote by $\mathcal{S}$ the set of all sectorial functions
(in $L^{\infty}(\mathbb{R})$). Once again, for Wiener-Hopf
operators with such kind of Fourier symbols a description of the
possible regularity properties is known. Additionally, more
detailed information (upon some factorizations of the original
Fourier symbols) is also known in the following form.

\begin{theorem}\cite[Theorem 6.17]{BoKaSp02} \label{thm4.2}
\label{theorem 1.2} If
$\varphi\in\mathcal{G}L^{\infty}(\mathbb{R})$, then:
\begin{itemize}
\item[(a)] $W_{\varphi}$ is invertible if and only if $\varphi= s h,
 s\in\mathcal{S},  h\in\mathcal{G}H^{\infty}_{\pm}\,;$\\[-3.003mm] \item[(b)]
$W_{\varphi}$ is left-invertible if and only if $\varphi= s h,\;
 s\in \mathcal{S},\;  h\in H^{\infty}_{+}\,;$\\[-3.003mm] \item[($b^{\prime}$)]
$W_{\varphi}$ is right-invertible if and only if $\varphi= s h,\;
 s\in \mathcal{S},\;  h\in H^{\infty}_{-}\,;$\\[-3.003mm] \item[(c)] $W_{\varphi}$
is Fredholm if and only if $\varphi= s h,\;  s\in \mathcal{S},\;
 h\in \mathcal{G}[C(\bkRp)+H^{\infty}_{\pm}]\,;$\\[-3.003mm] \item[(d)]
$W_{\varphi}$ is left-Fredholm if and only if $\varphi= s h,
 s\in\mathcal{S},\;  h\in C(\bkRp)+H^{\infty}_{+}\,;$\\[-3.003mm]
\item[($d^{\prime}$)] $W_{\varphi}$ is right-Fredholm if and only
if $\varphi= s h,  s\in\mathcal{S},\;  h\in
C(\bkRp)+H^{\infty}_{-}$.
\end{itemize}
\end{theorem}

Such type of theorem is also valid for Wiener-Hopf plus Hankel operators,
involving naturally some modifications over the symbol conditions. To prove the
next theorem in this direction, we first need to recall some auxiliary results.

\begin{theorem}[Pousson and Rabindranathan \cite{Po68, Ra69}] \label{theorem PR} If $a\in
\mathcal{G}L^{\infty}(\mathbb{R})$, then there exist a unitary-valued $u$ and
$ h\in\mathcal{G}H^{\infty}_{+}$ such that $a=u h$, almost everywhere on
$\mathbb{R}$.
\end{theorem}

\begin{lemma}\cite[Lemma 2.21]{BoSi90} \label{lemma 1.4}  If $E$ is a subset of
$L^{\infty}(\mathbb{R})$, $\varphi\in L^{\infty}(\mathbb{R})$ is
unitary-valued, and $\dist(\varphi,E)<1$, then there exist a
function $f\in E$ and a sectorial function $ s\in
\mathcal{G}L^{\infty}(\mathbb{R})$ such that $\varphi= s f$.
\end{lemma}

Now we introduce and prove a Wiener-Hopf plus/minus Hankel version
of the Theorem~\ref{theorem 1.2}.

\begin{theorem}\label{proposition 1.5}
Let $\varphi\in\mathcal{G}L^{\infty}(\mathbb{R})$.

\begin{itemize}
\item[(a)] $WH_{\varphi}$ and $W_{\varphi}-H_{\varphi}$ are both
invertible if and only if $$\varphi\widetilde{\varphi^{-1}}= s
h,\;
 s\in \mathcal{S},\;  h\in\mathcal{G}H_{\pm}^{\infty}\,.$$

\item[(b)] $WH_{\varphi}$ and $W_{\varphi}-H_{\varphi}$ are  both
left-invertible if and only if $$\varphi\widetilde{\varphi^{-1}}=
s h,\;  s\in\mathcal{S},\;  h\in H_{+}^{\infty}\,.$$

\item[(c)] $WH_{\varphi}$ and $W_{\varphi}-H_{\varphi}$ are both
right-invertible if and only if $$\varphi\widetilde{\varphi^{-1}}=
s h,\;  s\in\mathcal{S},\;  h\in H_{-}^{\infty}\,.$$

\item[(d)] $WH_{\varphi}$ and $W_{\varphi}-H_{\varphi}$ are both
Fredholm if and only if $$\varphi\widetilde{\varphi^{-1}}= s h_{1}
h_{2},\;
 s\in\mathcal{S},\;  h_{1}\in \mathcal{G}[C(\bkRp)+H_{\pm}^{\infty}],\; h_{2}\in\mathcal{G}H_{\pm}^{\infty}.$$
\item[(e)] $WH_{\varphi}$ and $W_{\varphi}-H_{\varphi}$ are both
left-Fredholm if and only if $$\varphi\widetilde{\varphi^{-1}}= s
h_{1} h_{2},\;  s\in\mathcal{S},\;  h_{1}\in
C(\bkRp)+H_{+}^{\infty},\,\, h_{2}\in\mathcal{G}H_{+}^{\infty}\,.$$

\item[(f)] $WH_{\varphi}$ and $W_{\varphi}-H_{\varphi}$ are both
right-Fredholm if and only if $$\varphi\widetilde{\varphi^{-1}}= s
h_{1} h_{2},\;
 s\in\mathcal{S},\;  h_{1}\in C(\bkRp)+H_{-}^{\infty},\;
h_{2}\in\mathcal{G}H_{-}^{\infty}\,.$$
\end{itemize}
\end{theorem}

\begin{proof}
First we will prove the ``+"-version (cf.~the indices). As a
global property, from the ${\Delta}$-relation after extension we
know that
$WH_{\varphi}\stackrel{*}{\Delta}W_{\varphi\widetilde{\varphi^{-1}}}$.
 Additionally, from the Theorem \ref{theorem PR} (of Pousson and
Rabindranathan), we have
\begin{eqnarray}\label{eqqa}
\varphi\widetilde{\varphi^{-1}}=u_{\varphi} h_{\varphi}\;,
\end{eqnarray}
where $u_{\varphi}$ is unitary-valued, and $ h_{\varphi}\in
\mathcal{G}H^{\infty}_{+}$ (the notation of $\varphi$ in index is only to
reinforce that those factors depend on $\varphi$). Therefore, (\ref{eqqa})
yields that $W_{\varphi\widetilde{\varphi^{-1}}}$ is equivalent to
$W_{u_{\varphi}}$. Using now the Douglas--Sarason Theorem
(cf.~Theorem~\ref{theorem1}), we have:
\begin{itemize}
\item[(i)] $W_{u_{\varphi}}$ is invertible if and only if
$\dist(u_{\varphi},\mathcal{G}H_{+}^{\infty})<1\,;$\\[-3.003mm]

 \item[(ii)] $W_{u_{\varphi}}$ is left-invertible if and only if
 $\dist(u_{\varphi},H^{\infty}_{+})<1\,;$\\[-3.003mm]

\item[(iii)] $W_{u_{\varphi}}$ is Fredholm if and only if
$\dist(u_{\varphi},\mathcal{G}[C(\bkRp)+H_{+}^{\infty}])<1\,;$\\[-3.003mm]

\item[(iv)] $W_{u_{\varphi}}$ is left-Fredholm if and only if
$\dist(u_{\varphi},C(\bkRp)+H_{+}^{\infty})<1\,.$\\[-0.5mm]
\end{itemize}
This combined with the above Lemma~\ref{lemma 1.4} yields that:
\begin{itemize}
\item[(i$^\prime$)] If $W_{u_{\varphi}}$ is invertible, then $u_{\varphi}= s h$
for some $s\in\mathcal{S}$ and $h\in\mathcal{G}H_{+}^{\infty}\,;$\\[-3.003mm]
\item[(ii$^{\prime}$)] If $W_{u_{\varphi}}$ is left-invertible, then
$u_{\varphi}= s h$ for some $s\in\mathcal{S}$ and $h\in
H_{+}^{\infty}\,;$\\[-3.003mm]
\item[(iii$^\prime$)] If $W_{u_{\varphi}}$ is Fredholm, then $u_{\varphi}= s h$
for some $s\in \mathcal{S}$ and
$h\in\mathcal{G}[C(\bkRp)+H_{+}^{\infty}]\,;$\\[-3.003mm]
\item[(iv$^\prime$)] If $W_{u_{\varphi}}$ is left-Fredholm, then
$u_{\varphi}= s h$ for some $s\in \mathcal{S}$ and $h\in
C(\bkRp)+H_{+}^{\infty}\,$.\\[-0.5mm]
\end{itemize}
Using the multiplication by $ h_{\varphi}$ in the last four
identities (and recalling the identity
$\varphi\widetilde{\varphi^{-1}}=u_{\varphi} h_{\varphi}$), we
conclude that:
\begin{itemize}
\item[(i$^{\prime\prime}$)] If $W_{u_{\varphi}}$ is invertible, then
$$\varphi\widetilde{\varphi^{-1}}= s h, \quad  s\in\mathcal{S}, \;
 h\in\mathcal{G}H_{+}^{\infty}\,;$$
\item[(ii$^{\prime\prime}$)] If $W_{u_{\varphi}}$ is left-invertible, then
$$\varphi\widetilde{\varphi^{-1}}= s h, \quad s\in\mathcal{S}, \;
 h\in H_{+}^{\infty}\,;$$
\item[(iii$^{\prime\prime}$)] If $W_{u_{\varphi}}$ is Fredholm, then
$$\varphi\widetilde{\varphi^{-1}}= s h_{1} h_{2}, \quad
 s\in\mathcal{S}, \;  h_{1}\in\mathcal{G}[C(\bkRp)+
H_{+}^{\infty}],\; h_{2}\in\mathcal{G}H_{+}^{\infty}\,;$$
\item[(iv$^{\prime\prime}$)] If $W_{u_{\varphi}}$ is left-Fredholm, then
$$\varphi\widetilde{\varphi^{-1}}= s h_{1} h_{2}, \quad
 s\in\mathcal{S}, \;  h_{1}\in C(\bkRp)+H_{+}^{\infty},\;
h_{2}\in\mathcal{G}H_{+}^{\infty}\,$$
\end{itemize}
where $ h_{1}$ and $ h_{2}$ represent $ h$ and $ h_{\varphi}$,
respectively.

Now we will proceed with the reverse implications. Let us then assume the
following four situations:
\begin{itemize}
\item[(j)] $\varphi\widetilde{\varphi^{-1}}= s h$, for some $s\in\mathcal{S}$
and $h\in\mathcal{G}H_{+}^{\infty}$; \item[(jj)]
$\varphi\widetilde{\varphi^{-1}}= s h$, for some $s\in\mathcal{S}$ and $h\in
H_{+}^{\infty}$; \item[(jjj)] $\varphi\widetilde{\varphi^{-1}}= s h_{1} h_{2}$,
for some $s\in\mathcal{S}$,  $h_{1}\in\mathcal{G}[C(\bkRp)+ H_{+}^{\infty}]$
and $h_{2}\in\mathcal{G}H_{+}^{\infty}$; \item[(jv)]
$\varphi\widetilde{\varphi^{-1}}= s h_{1} h_{2}$, for some $s\in\mathcal{S}$,
$h_{1}\in C(\bkRp)+ H_{+}^{\infty}$ and $h_{2}\in\mathcal{G}H_{+}^{\infty}$.
\end{itemize}
In the first two cases (j) and (jj) we directly obtain from
Theorem~\ref{thm4.2} that $W_{\varphi\widetilde{\varphi^{-1}}}$ is
invertible or left-invertible, respectively. Therefore, the above
mentioned $\Delta$-relation after extension yields the same
properties for the diagonal matrix operator
$\mbox{diag}[WH_{\varphi}, W_{\varphi}-H_{\varphi}]$, and the
result follows (for these two cases).

 In the cases (jjj) and (jv), formula (\ref{eqqa}) yields that
$$W_{u_{\varphi}}=W_{s}\, W_{h_{1}}\, W_{h_{2}h_{\varphi}^{-1}}\,.$$
Note that $W_{s}$ is an invertible operator because $ s$ is
sectorial, and $W_{h_{2}h_{\varphi}^{-1}}$ is also invertible due
to $h_{2}h_{\varphi}^{-1}\in\mathcal{G}H_{+}^{\infty}$ (having by
inverse operator the corresponding Wiener-Hopf operator with
Fourier symbol $h_2^{-1}h_{\varphi}$). Therefore,
$W_{u_{\varphi}}$ is equivalent to $W_{h_{1}}$. Thus, to conclude
the desired result by the use of the $\Delta$-relation after
extension, we only need to describe the
 analysis of $W_{h_1}$ in those corresponding cases:
\begin{itemize}
\item[(k)] If $h_1\in\mathcal{G}[C(\bkRp)+H_{+}^{\infty}]$, then
$W_{h_1}$ is Fredholm;\\[-3.003mm]
\item[(kk)] If $h_1\in
C(\bkRp)+H_{+}^{\infty}$, then $W_{h_1}$ is left-Fredholm.
\end{itemize}
In fact, propositions (k) and (kk) arise from the {\it Douglas
Theorem} (cf., e.g., \cite[Theorem 5.7]{BoKaSp02}), and by
observing that $ h_1\in\mathcal{G}L^{\infty}(\mathbb{R})$ (because
$\varphi\widetilde{\varphi^{-1}}\in\mathcal{G}L^{\infty}(\mathbb{R})$,
$h_{2}\in\mathcal{G}H_{+}^{\infty}$ and
$s\in\mathcal{G}L^{\infty}(\mathbb{R}))$.

We have therefore proved the ``+"--version of the theorem. Now, by
passage to adjoints, the remaining ``-"--versions are obtained.
\end{proof}

The ``only if" part of Theorem~\ref{proposition 1.5}, directly
yields the main conclusion of the present work:

\begin{corollary}\label{last:Cor}
Let $\varphi\in\mathcal{G}L^{\infty}(\mathbb{R})$.
\begin{itemize}
\item[(a)] If $\varphi\widetilde{\varphi^{-1}}= s h$, with $s\in
\mathcal{S}$ and $h\in\mathcal{G}H_{\pm}^{\infty},\;$ then
$WH_{\varphi}$ is invertible.\\[-3.003mm]
\item[(b)] If $\varphi\widetilde{\varphi^{-1}}= s h$, with
$s\in\mathcal{S}$ and $h\in
H_{+}^{\infty},\;$ then $WH_{\varphi}$ is left-invertible.\\[-3.003mm]
\item[(c)] If $\varphi\widetilde{\varphi^{-1}}=  s h$, with
$s\in\mathcal{S}$ and $h\in H_{-}^{\infty}$, then $WH_{\varphi}$
is
right-invertible.\\[-3.003mm]
\item[(d)] If $\varphi\widetilde{\varphi^{-1}}= s h_{1} h_{2}$,
with $s\in\mathcal{S}$, $h_{1}\in
\mathcal{G}[C(\bkRp)+H_{\pm}^{\infty}],$ and
$h_{2}\in\mathcal{G}H_{\pm}^{\infty},$
then $WH_{\varphi}$ is Fredholm.\\[-3.003mm]
\item[(e)] If $\varphi\widetilde{\varphi^{-1}}= s h_{1} h_{2}$,
with $s\in\mathcal{S}$, $h_{1}\in C(\bkRp)+H_{+}^{\infty}$ and
$h_{2}\in\mathcal{G}H_{+}^{\infty}$, then $WH_{\varphi}$ is
left-Fredholm.\\[-3.003mm]
\item[(f)] If $\varphi\widetilde{\varphi^{-1}}= s h_{1} h_{2}$,
with $s\in\mathcal{S}$, $h_{1}\in C(\bkRp)+H_{-}^{\infty}$ and
$h_{2}\in\mathcal{G}H_{-}^{\infty}$, then $WH_{\varphi}$ is
right-Fredholm.
\end{itemize}
\end{corollary}

We end with an example showing the applicability of the last
result. Let us consider the Wiener-Hopf plus Hankel operator
$$WH_{\varphi_{p}}: L_{+}^{2}(\mathbb{R}) \rightarrow
L^{2}(\mathbb{R}_{+})\,,$$ with the particular Fourier symbol
$$\varphi_{p}(x) = (2+\sin x)\;e^{i\alpha x}\;,\;\;\;\;\;
x\in\mathbb{R}\;,$$
where $\alpha\in\mathbb{R}$ is a given
parameter. Direct computations show that
$$\varphi_{p}(x)\;\widetilde{\varphi_{p}^{-1}}(x) = \frac{2+\sin x}{2-\sin x}\; e^{2i\alpha x}\;.$$
So, if we choose $s_{p}(x)={(2+\sin x)}/{(2-\sin x)}$ and
$h_{p}(x)=e^{2i\alpha x}$, we see that $s_{p}\in\mathcal{S}$. This
occurs because ${1}/{3}\leq s_{p} \leq 3$, and therefore the range
of $s_{p}$ is contained in the right half-plane (which boundary
passes through the origin). Moreover, depending whether
$\alpha\geq 0$ or $\alpha\leq 0,$ we have $h_{p}\in
H_{+}^{\infty}$ or $h_{p}\in H_{-}^{\infty},$ respectively.
Therefore, applying Corollary~\ref{last:Cor} to
$$\varphi_{p}\;\widetilde{\varphi_{p}^{-1}} = s_{p}\;h_{p}\;,$$
we conclude that:
\begin{itemize}
\item[(a)] if $\alpha=0,$ then $WH_{\varphi_{p}}$ is invertible;\\[-3.003mm]
\item[(b)] if $\alpha>0,$ then $WH_{\varphi_{p}}$ is
left-invertible;\\[-3.003mm] \item[(c)] if $\alpha<0,$ then $WH_{\varphi_{p}}$
is right-invertible. \end{itemize}

\vspace{3.00mm}

\noindent{\bf Acknowledgements}: The work was supported by {\em
Funda\c{c}{\~a}o para a Ci{\^e}ncia e a Tecnologia} through {\em Unidade de
Investiga\c{c}{\~a}o Matem{\'a}tica e Aplica\c{c}{\~o}es} of University of Aveiro,
Portugal.


\begin{thebibliography}{xx}

\bibitem{BaTo04}
{E.L. Basor and T. Ehrhardt},  Factorization theory for a class of
Toeplitz $+$ Hankel operators, {\it J. Oper. Theory}, {\bf
51}(2004), 411--433.

\bibitem{BoKaSp02}
{A. B{\"o}ttcher, Yu.I. Karlovich, I.M. Spitkovsky}, {\it Convolution
Operators and Factorization of Almost Periodic Matrix Functions},
Birkh{\"a}user Verlag, Basel (2002).

\bibitem{BoSi90}
A. B\"ottcher, B. Silbermann, {\it Analysis of Toeplitz
Operators}, Springer-Verlag, Berlin (1990).

\bibitem{CaSp98}
{L.P. Castro and F.-O. Speck}, Regularity properties and
generalized inverses of delta-related operators, {\it Z. Anal.
Anwendungen}, {\bf 17} (1998), 577--598.


\bibitem{CaSp05pm}
{L.P. Castro and F.-O. Speck}, Inversion of matrix convolution
type operators with symmetry, {\it Port. Math. (N.S.)}, {\bf 62}
(2005), 193--216.


\bibitem{CaSpTe04mn}
{L.P. Castro, F.-O. Speck and F.S. Teixeira}, A direct approach to
convolution type operators with symmetry, {\it Math. Nachr.}, {\bf
269/270} (2004), 73--85.

\bibitem{CaSpTe03jie}
{L.P. Castro, F.-O. Speck and F.S. Teixeira}, Explicit solution of
a Dirichlet-Neumann wedge diffraction problem with a strip, {\it
J. Integral Equations Appl.}, {\bf 15} (2003), 359--383.

\bibitem{CaSpTe04otaa}
{L.P. Castro, F.-O. Speck and F.S. Teixeira}, On a class of wedge
diffraction problems posted by Erhard Meister, {\it Oper. Theory
Adv. Appl.}, {\bf 147} (2004), 213--240.


\bibitem{CaSpTe06ieot}
{L.P. Castro, F.-O. Speck and F.S. Teixeira}, Mixed boundary value
problems for the Helmholtz equation in a quadrant, {\it Integral
Equations Operator Theory}, 44~pp., in press.

\bibitem{CoKrTe03}
{A.C. Concei\c{c}{\~a}o, V.G. Kravchenko and F.S. Teixeira}, Factorization
of some classes of matrix functions and the resolvent of a Hankel
operator, in: {\it Factorization, Singular Operators and Related
Problems}, Kluwer Acad. Publ., Dordrecht (2003), pp. 101--110.

\bibitem{DoSa}
R.G. Douglas and D. Sarason, Fredholm Toeplitz operators, {\it
Proc. Amer. Math. Soc.}, {\bf 26} (1970), 117--120.

\bibitem{Eh04}
{T. Ehrhardt}, Invertibility theory for Toeplitz plus Hankel
operators and singular integral operators with flip, {\it J.
Funct. Anal.}, {\bf 208} (2004), 64--106.

\bibitem{KaSa01}
{N. Karapetiants and S. Samko}, {\it Equations with Involutive
Operators}, Birkh{\"a}user, Boston (2001).

\bibitem{KrLeLiTe95}
{V.G. Kravchenko, A.B. Lebre, G.S. Litvinchuk and F.S. Teixeira},
Fredholm theory for a class of singular integral operators with
Carleman shift and unbounded coefficients, {\it Math. Nachr.},
{\bf 172} (1995), 199--210.

\bibitem{KrLi94}
{V.G. Kravchenko and G.S. Litvinchuk}, {\it Introduction to the
Theory of Singular Integral Operators with Shift},  Kluwer
Academic Publishers Group, Dordrecht (1994).

\bibitem{LeMeTe92}
{A.B. Lebre, E. Meister and F.S. Teixeira},  Some results on the
invertibility of Wiener-Hopf-Hankel operators, {\it Z. Anal.
Anwendungen}, {\bf 11} (1992), 57--76.


\bibitem{NoCa04}
{A.P. Nolasco and L.P. Castro}, Factorization of Wiener--Hopf plus
Hankel operators with $APW$ Fourier symbols, {\it Int. J. Pure
Appl. Math.}, {\bf 14} (2004), 537--550.

\bibitem{NoCa05}
{A.P. Nolasco and L.P. Castro}, A criterion for lateral
invertibility of matrix Wiener-Hopf plus Hankel operators with
good Hausdorff sets, 9 p., {\it to appear}.

\bibitem{NoCa05b}
{A.P. Nolasco and L.P. Castro}, A Duduchava-Saginashvili's type
theory for Wiener-Hopf plus Hankel operators, 16 p., {\it to
appear}.


\bibitem{Po68}
H.R. Pousson,  Systems of Toeplitz operators on $H\sp{2}$, {\it
Proc. Amer. Math. Soc.}, {\bf 19} (1968), 603--608.

\bibitem{Ra69}
M. Rabindranathan, On the inversion of Toeplitz operators, {\it J.
Math. Mech.}, {\bf 19} (1969/1970), 195--206.


\bibitem{RoSi90}
{S. Roch and B. Silbermann},  {\it Algebras of Convolution
Operators and their Image in the Calkin Algebra},  {\it Akademie
der Wissenschaften der DDR, Karl-Weierstrass-Institut f{\"u}r
Mathematik}, Berlin (1990).


\bibitem{Te89}
{F.S. Teixeira}, On a class of Hankel operators: Fredholm
properties and invertibility, {\it Integral Equations Operator
Theory}, {\bf 12} (1989), 592--613.
\end{thebibliography}
\end{document}